\documentclass[a4paper,12pt]{article}
\usepackage[latin1]{inputenc}
\usepackage[T1]{fontenc}
\usepackage{amsfonts}
\usepackage[francais,english]{babel}
\usepackage{amssymb}
\usepackage{makeidx}
\usepackage{vmargin}
\usepackage{amsmath}
\usepackage{multirow}
\usepackage{amsthm}
\usepackage{graphics}
\usepackage{epsfig}
\newtheorem{theorem}{Theorem}
\newtheorem{lemma}{Lemma}
\newtheorem{e-proposition}{Proposition}
\newtheorem{corollary}{Corollary}
\newtheorem{e-definition}{Definition\rm}
\newtheorem{remark}{\it Remark\/}

\newtheorem{theoreme}{Th\'eor\`eme}

\newtheorem{proposition}{Proposition}

\newenvironment{demo}{{\sc Proof : }}{\begin{flushright}$\Box$\end{flushright}}

\begin{document}
\title{The centipede is determined by its Laplacian spectrum}
\author{Romain Boulet \\ Institut de Math\'ematiques de Toulouse\\ Universit\'e de Toulouse et CNRS (UMR 5219) \\ boulet@univ-tlse2.fr}
\date{This paper was submitted on February 22nd, 2008 in \textit{Comptes Rendus Mathématiques de l'Académie des Sciences de Paris}, accepted after revision on May 24th, 2008, available online on June 20th, 2008.}
\maketitle

\selectlanguage{english}
\begin{abstract}

A centipede is a graph obtained by appending a pendant vertex to each vertex of degree $2$ of a path. In this paper we prove that the centipede is determined by its Laplacian spectrum.
{\it To cite this article: R.
Boulet, C. R. Acad. Sci. Paris, Ser. I 346 (2008).}
\end{abstract}
\vskip 0.5\baselineskip

\selectlanguage{francais}
\begin{abstract}
\begin{center}{\bf Le mille-pattes est d\'etermin\'e par le spectre du Laplacien. }\end{center}
Un mille-pattes est un graphe obtenu en attachant un sommet pendant \`a chaque sommet de degr\'e $2$ d'une cha\^ine. Dans cet article nous montrons qu'un mille-pattes est d\'etermin\'e par le spectre du Laplacien.
{\it Pour citer cet article~: R. Boulet, C. R. Acad. Sci.
Paris, Ser. I 346 (2008).}

\end{abstract}

 \selectlanguage{francais}
 \section*{Version fran\c{c}aise abr\'eg\'ee}
 Le Laplacien $L$ d'un graphe est la matrice $L=D-A$ o\`u $D$ est la matrice diagonale des degr\'es et $A$ est la matrice d'adjacence du graphe.  Le spectre du Laplacien donne des informations sur la structure du graphe, comme la connexit\'e (voir  \cite{RePEc:dgr:kubcen:200266,mohar} pour plus de d\'etails), ces informations sont souvent insuffisantes pour reconstruire le graphe \`a partir du spectre  et  la question \og Quels graphes sont d\'etermin\'es par leur spectre ?\fg\ \cite{RePEc:dgr:kubcen:200266} demeure un probl\`eme difficile.
En particulier il est connu \cite{these_newman} que presque aucun arbre n'est d\'etermin\'e par le spectre du Laplacien et seules quelques familles d'arbres d\'etermin\'es par leur spectre ont jusqu'alors \'et\'e d\'ecouvertes (citons par exemple \cite{starlike_laplacian} et \cite{Z_n}).

On appelle mille-pattes le graphe obtenu en attachant un sommet pendant \`a chaque sommet de degr\'e $2$ d'une cha\^ine (voir figure \ref{chenille}). 
Nous montrons dans cet article que le mille-pattes est d\'etermin\'e par le spectre du Laplacien, enrichissant ainsi les familles connues d'arbres d\'etermin\'es par le spectre du Laplacien.

On note $P_k$ la cha\^ine \`a $k$ sommets et $T$ le triangle. Deux graphes sont dits $A$-cospectraux (resp. $L$-cospectraux) s'ils ont m\^eme spectre pour la matrice d'adjacence (resp. le Laplacien).  Les valeurs propres de la matrice d'adjacence d'un graphe $G$ d'ordre $n$ sont not\'ees $\lambda_1(G)\geq \lambda_2(G)\geq ...\geq \lambda_n(G)$ et celles du Laplacien $\mu_1(G)\geq \mu_2(G)\geq ... \geq\mu_n(G)$.
Le graphe repr\'esentatif des ar\^etes $\mathcal{L}(G)$ de $G$ a pour sommets les ar\^etes de $G$ et deux sommets sont adjacents si et seulement si les ar\^etes correspondantes dans $G$ poss\`edent un sommet commun.

Les r\'esultats suivants sont connus et permettent d'obtenir des informations sur la structure du graphe \`a partir du spectre du Laplacien:

\begin{theoreme}\label{th_fr}\cite{mohar}
Soit $G=(V,E)$ un graphe et soit $d(v)$ le degr\'e d'un sommet $v$. Alors :
$$\max\{d(v), v\in V(G)\}<\mu_1(G)\leq\max\{d(u)+d(v), \ uv\in E(G)\}$$
\end{theoreme}

\begin{theoreme}\label{laplacien_pptes_fr}\cite{RePEc:dgr:kubcen:200266}
Le nombre de sommets, d'ar\^etes, de composantes connexes et d'arbres couvrants d'un graphe peuvent \^etre d\'eduits du spectre de son Laplacien.
\end{theoreme}

\begin{theoreme}\cite{Doob_topic,mohar}\label{line_fr}
Soit $G$ un arbre \`a $n$ sommets, alors $\mu_i(G)=\lambda_i(\mathcal{L}(G))+2$ pour $1\leq i\leq n-1$.
\end{theoreme}

Le spectre de la matrice d'adjacence donne \'egalement des informations sur la structure du graphe; il est en particulier connu que $\sum_i \lambda_i^k$ est \'egal au nombre de marches ferm\'ees de longueur $k$ dans $G$. 

Soit $M$ un graphe, une marche ferm\'ee couvrante de longueur $k$ sur $M$ est une marche ferm\'ee de longueur $k$ sur $M$ parcourant \emph{toutes} les ar\^etes de $M$ au moins une fois. On note $w_k(M)$ le nombre de marches ferm\'ees couvrantes de longueur $k$ sur  $M$ et on d\'efinit l'ensemble $\mathcal{M}_k=\{M, w_k(M)>0\}$. Le nombre sous graphes (non n\'ecessairement induits) de $G$ isomorphes \`a $M$  est not\'e $|M(G)|$.

Ainsi le nombre de marches ferm\'ees de longueur $k$ de $G$ est : 
$$\sum_{M\in\mathcal{M}_k} w_k(M)|M(G)|$$

Une premi\`ere \'etape dans la d\'emonstration consiste \`a \'etudier l'ensemble $\mathcal{T}$ des graphes d\'efinis ainsi : $G$ est le graphe $T$ ou $G$ est obtenu en identifiant un sommet de degr\'e $2$ de $H\in\mathcal{T}$ et un sommet de  $T$.

Le sous-ensemble $\mathcal{T}_n$ de $\mathcal{T}$ constitu\'e des graphes de $\mathcal{T}$ avec exactement $n$ triangles est en bijection avec l'ensemble $\mathcal{A}_n$ des arbres \`a $n$ sommets de degr\'e maximal inf\'erieur ou \'egal \'a $3$. Cette bijection est l'application $K:\mathcal{T}_n\to\mathcal{A}_n$ 
qui \`a un graphe $G$ associe son graphe des cliques $K(G)$. Les sommets de $K(G)$ sont les sous-graphes complets maximaux de $G$ et deux sommets de $K(G)$ sont adjacents si et seulement si l'intersection des sous graphes complets maximaux correspondants dans $G$ est non vide.

Nous montrons ensuite qu'aucun graphe de $\mathcal{T}$ ne peut \^etre $A$-cospectral  avec $K^{-1}(P_k)$. Pour cela nous d\'enombrons le nombre de marches ferm\'ees de longueur $7$ pour $G\in\mathcal{T}\setminus\{T\}$ et nous obtenons $$\sum_i \lambda_i^7=686t-672+112t_3$$ o\`u $t$ est le nombre de sommets de $K(G)$ et $t_3$ le nombre de sommets de degr\'e $3$ de $K(G)$. Il reste \`a remarquer que $K^{-1}(P_k)$ minimise cette quantit\'e car $K(K^{-1}(P_k))=P_k$ ne poss\`ede aucun sommet de degr\'e $3$.

La deuxi\`eme \'etape de la d\'emonstration est de d\'eterminer la distribution des degr\'es d'un graphe $G$ $L$-cospectral avec un mille-pattes, on note $n_i$ le nombre de sommets de $G$ de degr\'e $i$. Le th\'eor\`eme \ref{th_fr} implique que le degr\'e maximal de $G$ est inf\'erieur ou \'egal \`a $5$. Puis  le th\'eor\`eme \ref{laplacien_pptes_fr} et la proposition \ref{sum_fr}
\begin{proposition}\label{sum_fr}
i) Soit $G$ un graphe, la somme des carr\'es des degr\'es de $G$ peut \^etre d\'eduite du spectre du Laplacien.
\\ ii) Soit $G$ un graphe dont on conna\^it le nombre de triangles, alors la somme des cubes des degr\'es de $G$ peut \^etre d\'eduite du spectre du Laplacien.
\end{proposition}
nous permettent d'obtenir le syst\`eme suivant que l'on r\'esoud.
$$\left\{\begin{array}{c}
n_1+n_2+n_3+n_4+n_5=n\\
n_1+2n_2+3n_3+4n_4+5n_5=2n-2\\
n_1+4n_2+9n_3+16n_4+25n_5=5n-8\\
n_1+8n_2+27n_3+64n_4+125n_5=14n-26
\end{array}\right.$$
 On obtient alors que $G$ est un arbre d'ordre $n$ avec $\frac{n+2}{2}$ sommets de degr\'e $1$ et $\frac{n-2}{2}$ sommets de degr\'e $3$. Il en d\'ecoule que $\mathcal{L}(G)\in\mathcal{T}$. Or, comme $G$ est $L$-cospectral avec un mille-pattes, le th\'eor\`eme \ref{line_fr} implique que $\mathcal{L}(G)$ est $A$-cospectral avec $K^{-1}\left(P_{\frac{n-2}{2}}\right)$ et donc $\mathcal{L}(G)$ est isomorphe \`a $K^{-1}\left(P_{\frac{n-2}{2}}\right)$ (premi\`ere \'etape de la d\'emonstration). On conclut en remarquant qu'un arbre avec $\frac{n+2}{2}$ sommets de degr\'e $1$ et $\frac{n-2}{2}$ sommets de degr\'e $3$ dont le graphe des lignes est isomorphe \`a $K^{-1}\left(P_{\frac{n-2}{2}}\right)$ est n\'ecessairement un mille-pattes.

\selectlanguage{english}

\section{Introduction}
\label{}

The Laplacian of a graph is the matrix $L=D-A$ where $A$ is the adjacency matrix and $D$ is the diagonal matrix of degrees. Some structural properties of the graph such as connectivity can be determined from the Laplacian spectrum  (see \cite{RePEc:dgr:kubcen:200266,mohar} for more details). However the Laplacian spectrum does generally not determine the graph and the question \emph{"Which graphs are determined by their spectrum ?"} \cite{RePEc:dgr:kubcen:200266} remains a difficult problem.
Moreover it is known \cite{these_newman} that almost no trees are determined by their Laplacian spectrum and the few trees proved to be determined by their Laplacian spectrum (see for instance \cite{starlike_laplacian} or \cite{Z_n}) may be viewed in this sense as exceptionnal graphs.

\begin{figure}[htbp]
\begin{center}
\includegraphics[scale=0.45]{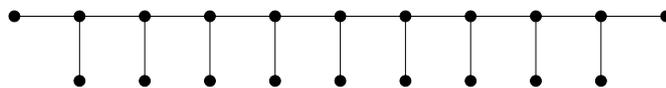}
\caption{A centipede}
\label{chenille}
\end{center}
\end{figure}

A centipede is a tree constructed by appending a pendant vertex to each vertex of degree $2$ of a path (see figure \ref{chenille} for an example). We show in this article that a centipede is determined by its Laplacian spectrum, thus enlarging the known families of trees determined by their Laplacian spectrum.
 
To fix notations, the path with $k$ vertices is denoted by $P_k$ and the triangle by $T$. Two graphs are $A$-cospectral (resp. $L$-cospectral) if they have the same adjacency (resp. Laplacian) spectrum. For a graph $G$ of order $n$, the eigenvalues of the adjacency matrix are denoted by $\lambda_1(G)\geq \lambda_2(G)\geq ...\geq \lambda_n(G)$ and the eigenvalues of the Laplacian by $\mu_1(G)\geq \mu_2(G)\geq ...\geq \mu_n(G)$ (when no confusion is possible we may omit to precise the graph).

We denote by $\textrm{Sp}(G)$ the spectrum of the adjacency matrix of $G$.
The line graph $\mathcal{L}(G)$ of a graph $G$ has the edges of $G$ as its vertices and two vertices of $\mathcal{L}(G)$ are adjacent if and only if the corresponding edges in $G$ have a commun vertex. For a vertex $v$ of a graph, $N(v)$ denotes the set of vertices adjacent to $v$ and $d(v)=|N(v)|$ the degree of $v$.

Here are some known results about the Laplacian spectrum.
\begin{theorem}\label{th}\cite{mohar}
Let $G=(V,E)$ be a graph where $V$ (resp. $E$) is the set of vertices (resp. edges). Then:
$$\max\{d(v), v\in V\}<\mu_1(G)\leq\max\{d(u)+d(v), \ uv\in E\}$$
\end{theorem}

\begin{theorem}\label{laplacien_pptes}\cite{RePEc:dgr:kubcen:200266}
For the Laplacian matrix, the following can be deduced from the spectrum:
\begin{itemize}
\item The number of vertices.
\item The number of edges.
\item The number of connected components.
\item The number of spanning trees.
\end{itemize}
\end{theorem}

\begin{theorem}\cite{Doob_topic,mohar}\label{line}
Let $G$ be a tree with $n$ vertices, then $\mu_i(G)=\lambda_i(\mathcal{L}(G))+2$ for $1\leq i\leq n-1$.
\end{theorem}

The spectrum of the adjacency matrix also gives some informations about structural properties of the graph, for instance it is a classical result that the number of closed walks of length $k\geq 2$ is $\sum_i\lambda_i^k$. 

We describe  here a method to count the number of closed walks of given length in a graph.

Let $M$ be a graph, a \emph{$k$-covering closed walk} in $M$ is a closed walk of length $k$ in $M$ running through \emph{all} the edges at least once. Let $G$ be a graph, $M(G)$ denotes the set of all distincts subgraphs (not necessarily induced) of $G$ isomorphic to $M$ and $|M(G)|$ is the number of elements of $M(G)$. 
The number of $k$-covering closed walks in  $M$ is denoted by $w_k(M)$ and we define the set $\mathcal{M}_k=\{M,\ w_k(M)>0\}$. 
As a consequence, the number of closed walks of length $k$ in $G$ is: 

 \begin{center}
\begin{equation}
\sum_{\lambda_i\in\textrm{Sp}(G)} \lambda_i^k=\sum_{M\in\mathcal{M}_k} w_k(M)|M(G)|
 \label{eq_closed_walk}
 \end{equation}
 \end{center}

\section{Preliminaries: definition and spectral properties of a set $\mathcal{T}$ of graphs }

Let $\mathcal{T}$ be the set of graphs $G$ defined as follow: $G$ is the triangle $T$ or $G$ is formed by identifying a vertex of degree $2$ of a graph $H\in\mathcal{T}$ and a vertex of the graph $T$.

The clique graph of $G$, denoted by $K(G)$, is the graph whose vertex set is the set of maximal complete subgraphs of $G$ and two vertices of $K(G)$ are adjacent if and only if the corresponding complete subgraphs in $G$ share at least one vertex.
The set of graphs in $\mathcal{T}$ with $n$ triangles is denoted by $\mathcal{T}_n$.

Let $\mathcal{A}_n$ be the set of trees on $n$ vertices with maximal degree lower than or equal to $3$. The following proposition is straightfoward:

\begin{e-proposition}
The application $K:\mathcal{T}_n\longrightarrow \mathcal{A}_n$ is one-to-one.
\end{e-proposition}

Equation (\ref{eq_closed_walk}) enables us to compute the number of closed walks of length $7$ of a graph $G\in\mathcal{T}$:

\begin{lemma}
We have for $G\in\mathcal{T}\setminus\{T\}$:
$$\sum_{\lambda_i\in\textrm{Sp}(G)} \lambda_i^7=686t-672+112t_3$$
where $t$ is the number of vertices of $K(G)$ (\textit{ie} the number of triangles in $G$)
 and $t_3$ is the number of vertices of degree  $3$ in $K(G)$.
\end{lemma}

\noindent\begin{demo}
We use the relation 
$$\sum_{\lambda_i\in\textrm{Sp}(G)} \lambda_i^7=\sum_{M\in\mathcal{M}_7} w_7(M)|M(G)|.$$

\begin{figure}[h]
\begin{center}
\includegraphics[scale=0.45]{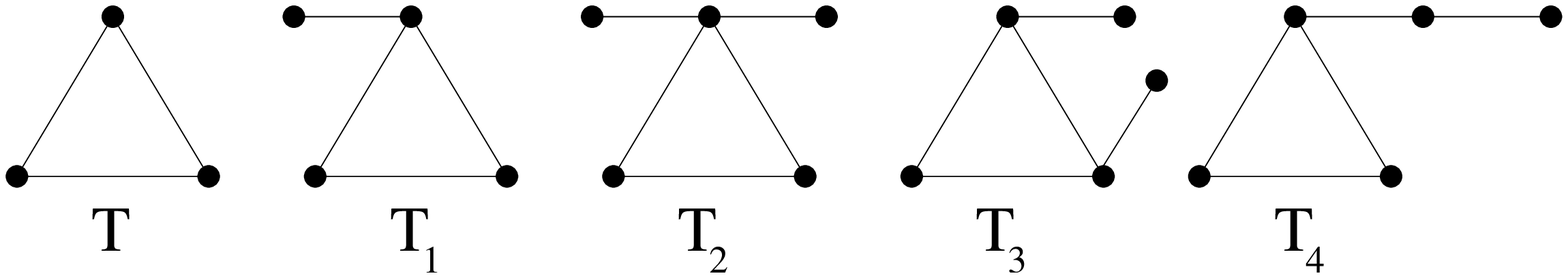}
\caption{subgraphs of $G\in\mathcal{T}$ belonging to $\mathcal{M}_7$}
\label{Ti}
\end{center}
\end{figure}

As an odd closed walk necessarily runs through an odd cycle, it is clear that $M\in\mathcal{M}_7$ contains one and only one triangle. Only the graphs  $T,\ T_1,\ T_2,\ T_3,\ T_4 \in \mathcal{M}_7 $ depicted in figure  \ref{Ti} can arise as subgraphs of $G\in \mathcal{T}$.

Let $t$ be the number of vertices of $K(G)$ and $t_i$, $1\leq i\leq3$, be the number of vertices of $K(G)$ of degree $i$. Since $K(G)\in\mathcal{A}_t$, we have $t_1+t_2+t_3=t$ and $t_1+2t_2+3t_3=2t-2$ (the sum of the degrees is twice the number of edges) showing $t_1=t_3+2$ and $t_2=t-2-2t_3$.

For a triangle $T$ of $G$ (denoted by $T\subset G$), let $N(T)$ be the set of triangles of $G$ sharing one vertex with $T$ and $d(T)=|N(T)|$. Note that $d(T)$ is the degree of the vertex in $K(G)$ corresponding to $T\subset G$. We have:
\begin{itemize}
\item $|T(G)|=t$.
\item $|T_1(G)|=2\sum_{T\subset G}d(T)=2\sum_{v\in V(K(G))}d(v)=2(2t-2)$
\item $|T_2(G)|=\sum_{T\subset G}d(T)=\sum_{v\in V(K(G))}d(v)=2t-2$
\item $|T_3(G)|=4t_2+12t_3=4(t-2)+4t_3$
\item $|T_4(G)|=4(2t-3+t_3)$, indeed  let $e$ and $e'$ be the two bridges of $T_4$ such that $e$ shares a vertex with the triangle of $T_4$.  Let $T$ be a triangle of $G$, then the number of $T_4$ such that $e$ and $e'$ belongs to $T$ is $2d(T)$. The number of $T_4$ such that  $e$ belongs to $T$ and $e'$ does not is $4$ if $d(T)=2$ or $12$ if $d(T)=3$. This implies $|T_4(G)|=(2t_1+4t_2+6t_3)+(4t_2+12t_3)$, that is, $|T_4(G)|=4(2t-3+t_3)$.
\end{itemize}

Let $A$, $A_i$ denote the adjacency matrices of $T$, $T_i$. The resolution of the equations $\textrm{tr}(A^7)=w_7(T)$, $\textrm{tr}(A_1^7)=w_7(T)+w_7(T_1)$, $\textrm{tr}(A_2^7)=w_7(T)+2w_7(T_1)+w_7(T_2)$, $\textrm{tr}(A_3^7)=w_7(T)+2w_7(T_1)+w_7(T_3)$, and $\textrm{tr}(A_4^7)=w_7(T)+w_7(T_1)+w_7(T_4)$ yields $w_7(T)=126$, $w_7(T_1)=84$, $w_7(T_2)=28$, $w_7(T_3)=14$ and $w_7(T_4)=14$.

This implies $\sum_i \lambda_i^7=\sum_{M\in\mathcal{M}_7} w_7(M)|M(G)|=686t-672+112t_3$.
\end{demo}

\begin{theorem}\label{iso}
The $A$-spectrum characterises $K^{-1}(P_k)$ in $\mathcal{T}$.
\end{theorem}

\noindent\begin{demo}
If $G\in\mathcal{T}$ is $A$-cospectral with $K^{-1}(P_k)$ then $G$ and $K^{-1}(P_k)$ have the same number of vertices and triangles. If $G$ and $K^{-1}(P_k)$ are not isomorphic then $K(G)$ possesses a vertex of degree $3$ (so $G$ is not isomorphic to $T$) and the previous lemma implies $\sum_{\lambda_i\in\textrm{Sp}(K^{-1}(P_k))} \lambda_i^7 < \sum_{\lambda_i\in\textrm{Sp}(G)} \lambda_i^7$. This  contradicts $A$-cospectrality of $G$ with $K^{-1}(P_k)$. 
\end{demo}
\section{The centipede is determined by its Laplacian spectrum}
\begin{proposition}
i) The sum of squares of the degrees of a graph $G$ can be deduced from its Laplacian spectrum.
\\ ii) The sum of cubes of the degrees of a graph $G$ can be deduced from its Laplacian spectrum and from the number of triangles contained in $G$.
\end{proposition}

\noindent\begin{demo}
i) We have 
$tr(L^2)=tr(D^2)-2tr(AD)+tr(A^2)$, but $tr(AD)=0$, $tr(A^2)=2m$  where $m$ is the number of edges 
and $tr(D^2)$ is the sum of squares of degrees of $G$.
\\
ii) We have 
$tr(L^3)=tr(D^3)-tr(A^3)+3tr(A^2D)$. But  $tr(A^3)$ is six times the number of triangles of $G$, $tr(A^2D)$ is the sum of squares of degrees of $G$ and $tr(D^3)$ is the sum of cubes of $G$.
\end{demo}

\begin{proposition}\label{distr_degres}
If $G$ is a graph on  $n$ vertices $L$-cospectral with a centipede,  then $G$ is a tree having $\frac{n-2}{2}$ vertices of degree $3$ and $\frac{n+2}{2}$ vertices of degree $1$.
\end{proposition}

\noindent\begin{demo}
Theorem \ref{th} implies that the maximal degree of $G$ is at most $5$. For $i=1,...,5$, let $n_i$ be the number of vertices of degree $i$ of $G$.
The Laplacian spectrum of $G$ determines the number of vertices, the number of edges, the sum of squares of the degrees and the sum of cubes of the degrees, that is:
$$\left\{\begin{array}{c}
n_1+n_2+n_3+n_4+n_5=n\\
n_1+2n_2+3n_3+4n_4+5n_5=2n-2\\
n_1+4n_2+9n_3+16n_4+25n_5=5n-8\\
n_1+8n_2+27n_3+64n_4+125n_5=14n-26
\end{array}\right.$$

This system implies $n_2=n_4=-4n_5=2(n+2-2n_1)$ showing $n_2=n_4=n_5=0$, $n_3=\frac{n-2}{2}$, $n_2=0$ and $n_1=\frac{n+2}{2}$.
\end{demo}

\begin{proposition}\label{line_chenille}
Let $G$ be a tree with $\frac{n-2}{2}$ vertices of degree $3$ and $\frac{n+2}{2}$ vertices of degree $1$, if $\mathcal{L}(G)$ is isomorphic to $K^{-1}(P_k)$ then $G$ is a centipede.
\end{proposition}

\noindent\begin{demo}
If $G$ is not a centipede then it contains the subgraph $H$ depicted in figure \ref{H}. This implies that $\mathcal{L}(G)$ is not isomorphic to $K^{-1}(P_k)$.
\begin{figure}[htbp]
\begin{center}
\includegraphics[scale=0.45]{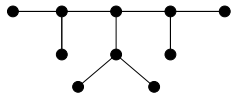}
\caption{Subgraph of a tree with $\frac{n-2}{2}$ vertices of degree $3$ and $\frac{n+2}{2}$ vertices of degree $1$ and different from a centipede }
\label{H}
\end{center}
\end{figure}
\end{demo}

\begin{theorem}
The centipede is determined by its Laplacian spectrum.
\end{theorem}

\begin{demo}
Let $G$ be a graph with $n$ vertices $L$-cospectral with the centipede on $n$ vertices. According to proposition  \ref{distr_degres}, $G$ is a tree with  $\frac{n-2}{2}$ vertices of degree $3$ and $\frac{n+2}{2}$ vertices of degree $1$. As $G$ is $L$-cospectral with a centipede, theorem \ref{line} implies that $\mathcal{L}(G)$ is $A$-cospectral with the line graph of the centipede, that is,  $\mathcal{L}(G)$ is $A$-cospectral with $K^{-1}\left(P_{\frac{n-2}{2}}\right)$.
Moreover $\mathcal{L}(G)\in\mathcal{T}$, so $\mathcal{L}(G)$ is isomorphic to  $K^{-1}\left(P_{\frac{n-2}{2}}\right)$ (theorem \ref{iso}) and  $G$ is a centipede by proposition \ref{line_chenille}.
\end{demo}

Since the Laplacian eigenvalues of a graph gives the Laplacian eigenvalues of its complement \cite{mohar},  we have the following corollary:
\begin{corollary}
The complement of a centipede is determined by its Laplacian spectrum.
\end{corollary}

\begin{remark}
 Non-uniquiness in $\mathcal{T}_n$ of graphs maximising $t_3$ prevents unfortunately the adaptation of our proof to graphs in $\mathcal{T}_n$ with $t_3$ maximal.
\end{remark}

\section*{Acknowledgements}
The author would like to thank the referee for his attentive reading and his relevant remarks.

\end{document}